\documentclass[runningheads]{amsart}
\usepackage{amssymb}
\usepackage{graphicx}
\usepackage{amsmath}
\usepackage{amscd}
\usepackage{enumerate}
\usepackage[none]{hyphenat}
\newif\ifskip
\skiptrue

\newif\ifshort
\shortfalse

\newif\ifcsl
\cslfalse

\newif\ifmargin
\marginfalse

\newtheorem{fact}{Fact}
\newtheorem{facts}{Facts}
\newtheorem{question}{Question}
\newtheorem{problem}{Problem}

\newtheorem{Conjecture}{Conjecture}

\newtheorem{theorem}{Theorem}[section]
\newtheorem{lemma}[theorem]{Lemma}
\newtheorem{proposition}[theorem]{Proposition}
\newtheorem{corollary}[theorem]{Corollary}
\newtheorem{definition}[theorem]{Definition}
\newtheorem{remark}[theorem]{Remark}
\newtheorem{example}[theorem]{Example}
\newtheorem{examples}[theorem]{Examples}

\newcommand{\N}{{\mathbb N}}
\newcommand{\Z}{{\mathbb Z}}
\newcommand{\R}{{\mathbb R}}
\newcommand{\Q}{{\mathbb Q}}
\newcommand{\C}{{\mathbb C}}

\newcommand{\cP}{{\mathcal P}}
\newcommand{\cQ}{{\mathcal Q}}

\newcommand{\cG}{{\mathcal G}}
\newcommand{\cD}{{\mathcal D}}
\newcommand{\cE}{{\mathcal E}}
\newcommand{\cH}{{\mathcal H}}

\newcommand{\rIND}{{\mathrm{IND}}}
\newcommand{\rDOM}{{\mathrm{DOM}}}

\ifskip
\else

\newcommand{\R}{{\mathbb R}}
\newcommand{\Z}{{\mathbb Z}}

\newcommand{\C}{{\mathbb C}}

\newcommand{\bl}{\begin{Lemm}}
\newcommand{\be}{\begin{equation}}
\newcommand{\el}{\end{Lemm}}
\newcommand{\ee}{\end{equation}}
\newcommand{\bt}{\begin{Theo}}
\newcommand{\et}{\end{Theo}}
\newcommand{\bp}{\begin{Prop}}
\newcommand{\ep}{\end{Prop}}

\newcommand{\bc}{\begin{Cor}}
\newcommand{\ec}{\end{Cor}}

\newenvironment{proof}{{\bf Proof:}\ }{\hfill Q.E.D.\newline}
\fi 

\renewcommand{\bar}{\overline}

\ifskip
\else
\newtheorem{theorem}{Theorem}[section]
\newtheorem{lemma}[theorem]{Lemma}

\theoremstyle{definition}
\newtheorem{definition}[theorem]{Definition}

\theoremstyle{remark}
\newtheorem{remark}[theorem]{Remark}
\fi

\numberwithin{equation}{section}

\newif\ifred
\redfalse
\newif\ifnewpage
\newpagefalse
\newif\ifmargin
\marginfalse

\usepackage[]{color}
\newcommand{\red}{\color{red}}

\newcommand{\black}{\color{black}}

\begin{document}

\title{Distinctive Power and Comparability \\ of Harary Polynomials}

\author{Johann A. Makowsky}
\address{Faculty of Computer Science, 
\newline Technion--Israel Institute of Technology, Haifa, Israel
}
\email{janos@cs.technion.ac.il}

\maketitle
Version: 
December 26, 2025
%
\begin{abstract}
Let $\cP$ be a graph property. A $\cP$-coloring with at most $k$ colors is a coloring of the vertices of a simple graph $G$
such that each color class induces a graph in $\cP$.
Harary polynomials are generalizations of the chromatic polynomial for simple graphs based on conditional colorings.
We denote by $\chi_{\cP}(G; k)$  the number of $\cP$-colorings of $G$ with at most $k$ colors.
$\chi_{\cP}(G; k)$ is a polynomial in $\Z[k]$.
A first paper studying Harary polynomials systematically was published in 2021 by O.Herscovici, J.A. Makowsky and V. Rakita.
\cite{HMR}. It studies under which conditions on $\cP$ is $\chi_{\cP}(G; k)$ definable in Monadic Second Order Logic 
and under which conditions is $\chi_{\cP}(G; k)$ a chromatic invariant.

Let $\cP, \cQ$ be two graph properties.
Two graphs $G, H$ are $\cP$-mates if \\
$\chi_{\cP}(G; k) = \chi_{\cP}(H; k)$.
$\chi_{\cQ}$ is at least as distinctive as $\chi_{\cP}$, 
$\chi_{\cP} \leq \chi_{\cQ}$,
if for all graphs $G, H$  we have that
$\chi_{\cQ}(G; k) = \chi_{\cQ}(H; k)$ implies 
$\chi_{\cP}(G; k) = \chi_{\cP}(H; k)$.
In this paper we study under which conditions on $\cP$ are there any (many) $\cP$-mates
and under which conditions on $\cP, \cQ$ is $\chi_{\cQ}$ is at least as distinctive as $\chi_{\cP}$.
\end{abstract}

\section{Introduction}
\label{se:intro}
\ifmargin
\marginpar{\red File:
\\ CM-intro.tex
\\
Last updated:
\\
20.12.25}
\else \fi 
\black

This paper is devoted to the systematic study of a large class of graph polynomials, the {\em Harary polynomials}.
They were previously studied in \cite{goodall2018complexity,makowsky2019logician,Rakita-PhD,HMR}. 

\subsection{Harary polynomials}
\begin{definition}
\label{def:harpol}
Let  $\cP$ be a graph property of simple graphs  and $G$ a graph. 
\begin{enumerate}[(i)]
\item
A {\em conditional $\cP$ coloring} short {\em $\cP$-coloring} with at most $k$ colors
is a function $c: V(G) \rightarrow [k]$ such that every color class $c^{-1}(i) = V_i$  
induces a graph 
\\$G[V_i] \in \cP$.
\item
Let $c_i^{\cP}(G)$ be the number of $\cP$-colorings of $G$ with exactly $i$ colors.
\item
A \em{conditional $\mathcal{P}$-partition} is a partition $\pi$ of $V(G)$ such that
$G[B]\in\mathcal{P}$ for each block $B$ of $\pi$. 
\item
Let $h_i^{\mathcal{P}}(G)$ be the number
of conditional $\mathcal{P}$-partitions of $G$ with exactly 
\\
$i$ blocks. 
\item
Let $\hat{\chi}_{\cP}(G;k)$ be the number of $\cP$-colorings of $G$ with at most $k$ colors.
\item
Let $\chi_{\cP}(G;k)$ be the number of $\cP$-partitions of $G$ with at most $k$ blocks.
\end{enumerate}
\end{definition}

$\cP$-colorings were first defined as a general concept by M. Harary in 1985, \cite{pr:Harary85}, and further
studied in \cite{brown1987generalized}.
The following may have been noticed before, 
see also \cite{biggs1974algebraic}, but was first intensively  studied in \cite{makowsky2006polynomial,ar:KotekMakowskyZilber11}.

\begin{theorem}
\label{th:two-harary}
\label{pr:basic}
For all simple graphs $G,H$ we have
\begin{enumerate}[(i)]
\item
$c_i^{\cP}(G) = (i)! \cdot h_i^{\mathcal{P}}(G)$.
\item
$\chi_{\cP}(G;k)= \sum_{i=0}^{n} h_i^{\cP}(G) k_{(i)}$ is a polynomial in $\Z[k]$.
\item
$\hat{\chi}_{\cP}(G;k)= \sum_{i=0}^{n} c_i^{\cP}(G) k_{(i)}$ is a polynomial in $\Z[k]$.
\item
If $G$ is a graph of order $n$ the polynomials $\chi_{\cP}(G;k)$ and $\hat{\chi}_{\cP}(G;k)$
are uniquely determined
by the graph parameters $h_i^P(G), i \leq n$ and $c_i^P(G), i \leq n$ respectively.
\item
$\chi_{\cP}(G;k) = \chi_{\cP}(H;k)$ iff $\hat{\chi}_{\cP}(G;k) = \hat{\chi}_{\cP}(H;k)$
\end{enumerate}
\end{theorem}
The proof is left to the reader.

We call $\chi_{\cP}(G;k)$ the {\em partition Harary polynomial of $G$} or just the  {\em Harary \\ polynomial of $G$},
and $\hat{\chi}_{\cP}(G;k)$ the {\em coloring Harary polynomial of $G$}.
In \cite{HMR} and in this paper the polynomial  $\chi_{\cP}(G;k)$ is used to define the Harary polynomial. 
In other papers which also deal with Harary polynomials,  $\hat{\chi}_{\cP}(G;k)$ is called th Harary polynomial.
However, we shall see in Theorem \ref{pr:harary-dp} below,
that they distinguish the same graphs. By Theorem \ref{th:two-harary}(i) the coefficients of  $\chi_{\cP}(G;k)$ are somehow simpler to manipulate because
of the absence of the factorial factor.

\subsection{Distinguishing power of graph polynomials}
The main theme in this paper
is the {\em distinguishing power} of Harary polynomials. 

The definitions below can be applied to any graph parameter or graph polynomial, not only to Harary polynomials.
\begin{definition}[Distinctive power]
Let $F, F_1, F_2$ be graph parameters. 
\begin{enumerate}[(i)]
\item
Two non isomorphic graphs 
$G$ and $H$ are 
{\em $F$-mates} or alternatively are called $F$-equivalent, if $F(G)=F(H)$, and 
a graph $G$ is {\em $F$-unique} if $G$ has no $F$-mates. 
\item
$U_F(n)$ is the set of 
$F$ unique graphs of order $n$. 
\item
$U_F = \bigcup_n U_F(n)$ is the set of $F$-unique graphs.
\item
$F$ is {\em trivial} if all pairs of graphs $G,H$ are $F$-mates.
\item
$F$ is {\em complete} if all graphs $G$ are $F$-unique.
\item
$F_1$ is less distinguishing than $F_2$, denoted by $F_1 <_{d.p.} F_2$ 
if 
for all graphs $G,H$ we have $F_2(G) = F_2(H)$ implies $F_1(G) = F_1(H)$.
d.p. stands for {\em distinctive power}.
\item
$F_1$ and $F_2$ are {\em comparable}, if either $F_1 <_{d.p.} F_2$ or $F_2 <_{d.p.} F_1$, but not both.
\item
$F_1$ is {\em equally distinguishing as $F_2$}, denoted by $F_1 \sim_{d.p.} F_2$ 
if for all $n \in \N^+$ we have  $U_{F_1}(n) =  U_{F_2}(n)$.
We also say that $F_1$ and $F_2$ are d.p.-equivalent.
\item
Let $F$ be a graph parameter and $G$ a graph. We denote by $[G]_F$  the class of $F$-mates of $G$.
\end{enumerate}
\end{definition}
In \cite{liu1997adjoint,noy2003graphs} $F$-unique graphs are studied for the characteristic, matching,
chromatic and Tutte polynomials.

\begin{proposition}
\label{pr:dp}
Let $F_1, F_2$ be two graph parameters.
The following are equivalent:
\begin{enumerate}
\item
$F_1 <_{d.p.} F_2$
\item
For every two graphs $G, H$ we have \\ $F_2(G) = F_2(H)$  implies $F_1(G) = F_1(H)$. 
\item
If $G,H$ are $F_2$-mates, they are also $F_1$-mates.
\\
In other words $[G]_{F_2} \subseteq [G]_{F_1}$.
\end{enumerate}
In particular:
$U_{F_1}(n) \subseteq U_{F_2}(n)$
\end{proposition}
\begin{proof}
This follows from the definitions.
\end{proof}

\begin{proposition}
\label{pr:harary-dp}
The two versions of Harary polynomials $\hat{\chi}_{\cP}(G;k)$ $\chi_{\cP}(G;k)$ 
are d.p.-equivalent.
\end{proposition}
\begin{proof}
This follows from the definition of d.p.-equivalence and Theorem \ref{pr:basic}(v).
\end{proof}

\subsection{Goal of the paper and main results}
In this paper we discuss some open questions 
about the distinctive power of Harary polynomials. 
In some cases we even dare to formulate conjectures.
In order to illustrate this, we present also various partial results.
The paper is an illustration  of how to the develop a general more abstract theory of graph polynomials as
summarized in \cite{makowsky2025meta}.
In Section \ref{se:chromatic} the main results from \cite{HMR} are summarized. 
The main new results in this paper are
Theorem \ref{th:hereditary},
Theorem \ref{th:gen-hereditary},
Theorem \ref{th:dp-1},
Theorem \ref{th:dp-2} and
Theorem \ref{th:delta}.

\subsection{Outline of the paper}

In Section \ref{se:background} we give more background and definitions and collect simple facts about Harary polynomials.
This is mostly taken from \cite{HMR}.

In Section \ref{se:chromatic} we discuss the chromatic polynomial for simple loop-less graphs
from our point of view in order to give motivation for our further questions about Harary polynomials
and list also the main open questions and conjectures about them.

In Section \ref{se:complete} we show some evidence for Conjecture \ref{conj-1}
which states that no Harary polynomial is complete.

In Section \ref{se:hereditary} we study Conjecture \ref{conj-1} for Harary polynomials $\chi_{\cP}$ with 
$\cP$ hereditary.

In Section \ref{se:generating} we digress and show how our methods can also be applied to graph polynomials
which are not Harary polynomials, but rather generating functions for counting induced subgraphs with
a graph property $\cP$.

In Section \ref{se:compare} we study d.p.-incomparability of Harary polynomial.

In Section \ref{se:wd} we explore under what conditions almost all graphs of
a Harary polynomial $\chi_{\cP}$ are not $\chi_{\cP}$-unique.

Finally, in Section \ref{se:conclu} we draw our conclusions.
\ifnewpage
\newpage
\else \fi 
\section{Background}
\label{se:background}
\ifmargin
\marginpar{\red File: 
\\
CM-background.tex
\\
Last updated:
\\
16.12.25}
\else \fi 
\black

\subsection{Graph parameters}
Recall that $[n] = \{1,2, \dots, n\}$.
We denote by $\cG$ the class of all simple finite graphs $(V(G), E(G))$. 
To make it into a set we assume that  $V(G)= [n]$ for some $n \in \N^+$.
$\cG(n)$ is the family of graphs of order $n$ with $V(G)= [n]$. 

A {\em graph parameter $F$} is a function mapping finite graphs $\cG$ into a domain 
$D$ which is invariant under graph isomorphism.
If $D =\{0,1\}$ we call $F$ a {\em graph property} and denote it usually by $\cP, \cQ$ possibly 
indexed by  some suggestive letter. 
If $D$ is a set of numbers like $\N, \Z, \Q, \R$ or $\C$ we speak of $F$ as a {\em numeric graph parameters}.
If $D$ is a polynomial ring in a fixed finite number of indeterminates $\bar{x} = (x_1, \ldots, x_k)$ 
like $\N[\bar{x}], \Z[\bar{x}], \Q[\bar{x}], \R[\bar{x}]$ or $\C[\bar{x}]$  we speak of $F$ as a {\em graph polynomials}.
If $F$ is a graph polynomial with $\bar{x}=x$ we speak of {\em univariate} graph polynomials.
Typical examples are:
\begin{description}
\item[Graph properties] The class of connected graphs,  planar graphs, edge-less graphs, regular graphs, etc.
\item[Numeric graph parameters] The number of vertices $n(G)$, the number of edges $m(G)$,
the number of connected components $k(G)$, 
the maximal degree, 
the chromatic number $\chi(G)$, etc.
\item[Univariate graph polynomials] The chromatic polynomial $\chi(G;k)$, the matching polynomial, 
the independence polynomial $\rIND(G;x)$, \\
the characteristic polynomial $char(G;x)$  
are all univariate graph polynomials.
\item[Multivariate graph polynomials] The Tutte polynomial $T(G; x,y)$, the edge elimination polynomial $\xi(G; x,y,z)$
\end{description}
A monograph with a comprehensive study of graph polynomials by P. Tittmann is in preparation \cite{bk:Tittmann-book}
and should appear soon. A preliminary version is available as \cite{Tittmann-draft}.

\subsection{Examples of Harary polynomials}
Let  $\cP$ be a graph property and $G$ a graph. 
Recall the definition of the Harary polynomials 
$$\chi_P(G;k) = \sum_i^n h_i^P(G) k_{(i)},$$
where
$h_i^P(G)$ be the number $P$-colorings of $G$ with {\em exactly $i$ colors}.

The following are Harary polynomials for $\cP$-colorings previously studied in the literature.
\begin{examples}
\label{ex:harary}
\ 
\begin{description}
\label{chi-mates}
\item[$\cP =\cE$]
If $\cP = \cE$ where $\cE$ denotes the set of edge-less graphs, 
$\chi_{\cE}(G;k)$ is the classical chromatic polynomial $\chi(G;k)$.
The graphs $C_n$ (cycles of order $n$) and $K_n$ (complete graphs of order $n$) are $\chi$-unique.
All trees of order $n$ are $\chi$-mates.
\item[$\cP =\cG$]
If $\cP = \cG$ is the set of all finite graphs and $G$ is of order $n$, we have 
\begin{gather}
h_i^{\cG}(i) = S(n,i) 
\notag\\
\chi_{\cG}(G;k) = \sum_{i=0}^k S(n,i)k_{(i)} 
\notag\\
\chi_{\cG}(G;k) = n^k
\end{gather}
All graphs of order $n \geq 2$ are $\chi_{\cG}$-mates. 
\notag
$K_1$ is the only $\chi_{\cG}$-unique graphs of order $1$ without loops.
\item[$mcc_t$-colorings] 
Here $\cP = \cP_t$, where $\cP_t$ is the graph property such that $H \in \cP_t$ iff
the connected components of $H \in \cP_t$ have order at most $t$.
\\
For $t=1$ these are the proper colorings, for $t=2$ these are the $\cP_3$-free colorings.
They were introduced in \cite{linial2007graph} with a slightly different notation.
\item[$\cP = DU(H)$]
Let $H$ be a connected graph of order $t$.
$DU(H)$ consist of non-empty disjoint unions of copies of $H$.
They are studied in
\cite{goodall2018complexity}.
\item[$\cP = Fr(H)$]
Let $Fr(H)$ be the class of all graphs which do not contain $H$ as an induced subgraphs. These are the $H$-free graphs.
$Fr(H)$-colorings are studied in \cite{achlioptas1997complexity,brown1996complexity,brown1987generalized}.
$K_2$-free colorings are just the proper colorings.
\item[Adjoint polynomial]
The {\em adjoint polynomial $A(G;x) = \chi_{\cP}(G; k)$} is defined by taking  $\cP$ to be the 
class of all complete graphs.  It was introduced in \cite{liu1997adjoint}, see also \cite{bencs2017one}.
\item[$\cP$ additive and hereditary]
A graph property $\cP$ is additive if it is closed under forming disjoint unions of graphs.
$\cP$ is hereditary, if it is closed under taking induced subgraphs.
A coloring is $\mathcal{AH}$ if it is a $\cP$-coloring for some $\cP$ which is both additive and hereditary.
$\mathcal{AH}$-colorings were studied in \cite{farrugia2004vertex}.
\item[Compton-Gessel classes]
A graph property is a {Compton-Gessel class} if it is additive and closed under taking connected components.
They were introduced in \cite{Compton-PhD,compton1983some,gessel1984combinatorial}. 
$\mathcal{AH}$-classes are special cases of Compton-Gessel classes.
\end{description}
\end{examples}.

Here are a few simple facts shared by all Harary polynomials.
\begin{facts}
\label{facts}
For every graph property $\cP$ and every graph $G$ of order $n$ we have
\begin{enumerate}[(i)]
\item
$h_0^{\cP}(G)=0$, if the nullgraph $(\emptyset, \emptyset)$ is not in $\cP$.
\item
$h_1^{\cP}(G) \in \{0,1\}$ and $h_1^{\cP}(G) =1$ iff $G \in {\cP}$.
\item
$h_n^{\cP}(G) \in \{0,1\}$ and $h_n^{\cP}(G) =1$ iff $K_1 \in \cP$.
\\
In other words,
the polynomial $\chi_{\cP}(G; k)$ is monic of degree $|V(G)|$ iff $K_1 \in {\cP}$.
\item
If $k < n$ and $\chi_{\cP}(G; k) =0$ then for all $0 <\ell < k$ also $\chi_{\cP}(G, \ell) =0$.
\end{enumerate}
\end{facts}

\subsection{Graph polynomials which are not Harary polynomials}

The following is taken verbatim from \cite{HMR}.
Many familiar graph polynomials are not Harary polynomials of the form $\chi_{\cP}(G;x)$.
We generalize here Theorem 5.7 from  {\cite{makowsky2019logician}}.
\begin{lemma}
\label{le:1}
For every graph property $\cP$ we have
$$
\chi_{\cP}(G;1) =
\begin{cases}
1 & G \in \cP, \\
0 & G \not \in \cP.
\end{cases}
$$
\end{lemma}
Using Lemma \ref{le:1} we get
\begin{proposition}
\label{pr:notharary}
Let $F(G;x)$ be a graph polynomial and $G$ be a graph such that $F(G;1) \neq 0$ and $F(G;1) \neq 1$.
Then there is no graph property $\cP$ such that $\chi_{\cP}(G;x) = F(G;x)$.
\end{proposition}

The characteristic polynomial $char(G;x)$ of $G$ is the characteristic polynomial of its
adjacency matrix, and the Laplacian polynomial $Lap(G;x)$ is the characteristic polynomial of the
Laplacian matrix $G$, see \cite{bk:BrouwerHaemers2012}.

The matching polynomials are defined using $m_i(G)$, the number of matchings of $G$ of size $i$.
$$
M(G;x) = \sum_i m_i(G) x^i \mbox{ and } \mu(G;x) = \sum_i (-1)^i m_i(G) x^{n-2i}.
$$
$M(G;x)$ is the {\em generating matching polynomial} and $\mu(G;x)$ is the \\
{\em matching defect polynomial},
see 
\cite{lovasz2009matching}.

Let $\mathrm{in}_i(G)$ be the number of of independent sets of $G$ of size $i$, and
$\mathrm{d}_i(G)$ the number of dominating sets of $G$ of size $i$.
We define the {\em independence polynomial} $\rIND(G;x)$, \cite{pr:LevitMandrescu05},
and the {\em domination polynomial} $\rDOM(G;x)$,
\cite{arocha2000mean,phd:Alikhani,ar:KotekPreenSimonTittmanTrinks2012} as
$$
\rIND(G;x) = \sum_i \mathrm{in}_i(G) x^i  \mbox{ and }  \rDOM(G;x) = \sum_i \mathrm{d}_i(G) x^i.
$$ 
\begin{theorem}
The following are not Harary polynomials of the form $\chi_{\cP}(G;x)$:
\begin{enumerate}[(i)]
\item
The characteristic polynomial $\mathrm{char}(G;x)$ and the Laplacian $\mathrm{Lap}(G;x)$.
\item
The generating matching polynomial $M(G;x)$  and the defect matching \\
polynomial $\mu(G;x)$.
\item
The independence polynomial $\rIND(G;x)$.
\item
The domination polynomial $\rDOM(G;x)$.
\end{enumerate}
\end{theorem}
\begin{proof}
We use Proposition \ref{pr:notharary}.
\\
(i):
$char(C_4;x) = (x-2)x^2(x+2)$ and  $Lap(C_4;x) =  x(x-4)(x-2)^2$,
\\
hence $char(C_4;1) = Lap(C_4;1) =   -3$.

(ii):
$ M(C_4;x) =1 +4x +2x^2 \mbox{  and  }  \mu(C_4;x) = 1 +4x + 2x^2$,  
\\
hence $ M(C_4;1) =7 \mbox{  and  }  \mu(C_4;1) = 7$, 

(iii):
$ \rIND(C_4;x) = 1+ 4x +2x^2$, hence $\rIND(C_4;1) = 7$.

(iv):
$\rDOM(K_2;x) =2x+x^2$, hence  $\rDOM(K_2;1)= 3$.
\end{proof}

$\rDOM$ and $\rIND$ are special cases graph polynomials of the form
$$
\cP_{\Phi}(G;x) = \sum_{A \subseteq V(G): \Phi(A)} x^{|A|}.
$$
Graph polynomials of this form are {\em generating functions} 
counting subsets \\
$A \subseteq V(G)$ satisfying a property $\Phi(A)$.
In the cases above,
that $A$ is an independent, respectively a dominating set, see also \cite{makowsky2019logician}.
We say that {\em $\Phi$ determines $A$}, 
if for every graph $G$ there is a unique $A \subseteq V(G)$ which satisfies $\Phi(A)$.

\begin{theorem}
Assume that $\Phi$ does not determine $A$, then there is no graph property $\cP$ such that
for all graphs $G$
$\chi_{\cP}(G;x) = \cP_{\Phi}(G;x)$.
Hence $\cP_{\Phi}(G;x)$ cannot be a Harary polynomial.
\end{theorem}
\begin{proof}
By Lemma \ref{le:1} $\chi_{\cP}(G;1) \in \{0,1\}$ for all graphs $G$.
However, since $\Phi$ does not determine $A$,  there is a graph $H$ with $\cP_{\Phi}(H;1) \geq 2$.
\end{proof}

\ifnewpage
\newpage
\else \fi 

\section{Questions and conjectures}
\label{se:chromatic}
\ifmargin
\marginpar{\red File: 
\\
CM-chromatic.tex
\\
Last updated:
\\
17.12.25}
\else \fi 
We look first at some questions around the chromatic polynomial and then  formulate our conjectures.

\subsection{The chromatic polynomial}
\ \\
The following features of the chromatic polynomial motivated our study of Harary polynomials.
\begin{itemize}
\item
The chromatic polynomial of a graph $G$ is a chromatic invariant (see \ref{se:eei}),
i.e. satisfies a recurrence relation based on deletion and contraction of edges of $G$
\item
The chromatic polynomial has both infinitely $\chi$-mates and $\chi$-unique graphs, 
hence it is not complete (see Section \ref{se:complete} and \cite{bk:DongKohTeo2005}).
\item
It still might be almost complete 
as suggested in \cite{ar:BollobasPebodyRiordan2000}.
\item
By Proposition \ref{pr:dp} we have:
If two Harary polynomials $\chi_{\cP}$ and $\chi_{\cQ}$ have the same distinguishing power, they have the same unique graphs.
\end{itemize}

\begin{question}
Is there a graph property $\cP$ different from $\cE$ such that  $\chi_{\cP}$ and the chromatic polynomial have the same distinguishing power?
\end{question}

Here are two more questions about the chromatic polynomial from the literature 
which can be naively reformulated for arbitrary Harary polynomials.
However, we do not pursue these questions further in this paper.

\begin{itemize}
\item
The coefficients of the  chromatic polynomial  of a graph are log-concave, 
in particular their absolute values form a unimodal sequence.
This was conjectured by R.C. Read  in \cite{read1968introduction} for the chromatic polynomial of graphs and by
A. Heron, J.C. Rota and D. Welsh for the characteristic polynomial for matroids.
It was finally proved by J. Huh \cite{huh2012milnor,huh2012log}. For J. Huh's work and the conjecture on matroids
 the reader may consult G. Kalai 
\cite{kalai2022work}.
\begin{question}
For which graph properties $\cP$ do the coefficients of $\chi_{\cP}$ always form a log-concave sequence?
\end{question}
\begin{remark}
This might be a ``stupid'' question. The late Bernhard Neumann\footnote{
https://en.wikipedia.org/wiki/Bernhard\_Neumann
}
once accused me of asking a stupid question\footnote{
Is there a finitely presentable infinite group with only finitely many conjugacy classes, \cite{makowsky1974some}. 
The problem is still wide open and listed in \cite{baumslag2002open}.}.
A ``stupid'' question is, he said, a question where one has absolutely no clue how
to look for an answer. 
\end{remark}
\item
Herbert Wilf asked in 1973, \cite{wilf1973polynomials}:
Given a polynomial $F(x)$ of degree $n$ is there a graph $G$ of order $n$ such that $\chi(G; x) = F(x)$?
One can answer this question by computing by brute force $\chi(G;x)$ for all graphs with vertex set $V(G) =[n]$.
What Wilf may have had in mind was some list of algebraic properties of $F(x)$, or simply a feasible algorithm.
\begin{question}
Find criteria on $\cP$ for the existence (or non-existence)  of a feasible (polynomial time) algorithm which,
given a polynomial $F(x)$ of degree $n$, finds a graph $G$ of order $n$ such that $\chi(G; x) = F(x)$?
\end{question}
\end{itemize}

\subsection{Are there Harary polynomials which are chromatic invariants?}
\label{se:eei}

Let $G$ be a graph and $e$ an edge of $G$.
As usual we denote by $G_{-e}$ the graph obtained from $G$ by deleting $e$,
and by $G_{/e}$ the graph obtained from $G$ by contracting $e$ and replacing multiple edges by a single edge.
Following \cite[Chapter 9.1]{bk:Aigner2007},
a function $f$, which maps graphs into a polynomial ring $\mathcal{R}=\mathcal{F}[\bar{X}]$ with coefficients
in a field  $\mathcal{F}$ of characteristic $0$, is called a {\em chromatic invariant} (aka {\em Tutte-Grothendieck invariant})
if the following hold.

\begin{enumerate}[(i)]
\item
If $G$ has no edges, $f(G)=c_{|V(G)|}$.
\item
If  $e \in E(G)$ is a bridge, then $f(G)= A \cdot f(G_{-e})$. 
\item
If  $e \in E(G)$ is a loop, then $f(G)= B \cdot f(G_{-e})$. 
\item
There exist $\alpha, \beta \in \mathcal{R}$ such that
for every $e \in E(G)$, which is neither a loop nor a bridge, we have
$
f(G) = \alpha \cdot f(G_{-e}) + \beta \cdot f(G_{/e})
$.
\item
Multiplicativity: 
If $G = G_1 \sqcup G_2$ is the disjoint union of two graphs $G_1, G_2$ then
$
f(G) = f(G_1) \cdot f(G_2)
$.
\end{enumerate}
Here $c, A, B, \alpha, \beta$ are constants in $\mathcal{R}$.

\begin{theorem}[G. Birkhoff, 1910]
The chromatic polynomial $\chi(G;x)$ is a chromatic invariant with $\alpha =1, \beta = -1$ and $c =x$.
\end{theorem}

Chromatic invariants have a characterization using the Tutte polynomial \\
$T(G;x,y)$, see
\cite[Chapter 9.1, Theorem 9.5]{bk:Aigner2007}.

\begin{theorem}
Let $f$ be a chromatic invariant with $A, B, \alpha, \beta$ indeterminates as above. Then for all graphs
$G$ 
$$
f(G) =  
\alpha^{|E| - |V| + k(G)} \cdot \beta^{|V| -k(G)} \cdot T(G; \frac{A}{\beta}, \frac{B}{\alpha}).
$$
\end{theorem}
It follows by a counting argument that not all Harary polynomials are chromatic invariants.
In \cite{HMR} the following was shown:
\begin{theorem}[{\cite[Theorem 3.3]{HMR}}]
\label{th:noteei}
Let $\cP$ be a non-trivial hereditary (monotone, minor closed) graph property. 
Then $\chi_\cP$ is an chromatic invariant if and only if $\chi_\cP$ is the chromatic polynomial.
\end{theorem}

We consider $\cP$ being minor closed/monotone/hereditary as a globally defined feature of $\cP$. 
defined by forbidden minors, forbidden (induced) substructures,

Other globally defined features
are being definable in some fragment of Second Order Logic, 
being of bounded width (tree-width, clique-width, twin-width, etc..).
Harary polynomials definable in Monadic Second Order Logic are also discussed in \cite[Section 5]{HMR}.

\begin{question}
Is Theorem \ref{th:noteei} also true for other globally defined families of graph properties?
\end{question}

\subsection{Are there complete Harary polynomials?}
We have already mentioned that the chromatic polynomial is not complete: 
There are many pairs of graphs $G,H$ which are $\chi$-mates.
However, there are also infinitely many graphs $G$ which are $\chi$-unique.

Our first question therefore is whether there is a graph property $\cP$ such that $\chi_{\cP}$ is complete.
We conjecture that this is not the case. More precisely:

\begin{Conjecture}
\label{conj-1}
\ 
\begin{description}
\item[Weak form]
For every non-trivial graph property $\cP$ there exist at least two graphs $G, H$ such that
$\chi_{\cP}(G; k) = \chi_{\cP}(H; k)$.
In other words, $G$ and $H$ are $\chi_{\cP}$-mates.
\item[Strong form]
For every non-trivial $\cP$ there exist infinitely many pairs of graphs   $G_i, H_i, i \in \N$ such that
$\chi_{\cP}(G_i; k) = \chi_{\cP}(H_i; k)$.
In other words, $G_i$ and $H_i$ are $\chi_{\cP}$-mates for every $i \in \N$.
\end{description}
\end{Conjecture}

\subsection{Are there almost complete Harary polynomials?}
Bollob\'as, Pebody and Riordan, \cite{ar:BollobasPebodyRiordan2000}, asked whether the chromatic polynomial could be {\em almost complete}
in the sense that almost all graphs are $\chi$-unique.
To make this precise we have to specify what we mean by {\em almost all graphs}.
The simplest customary model for random graphs is $\mathbb{G}(n, p(n))$ where on graphs of order $n$ edges are chose with equal probability $p(n)$
where $p(n)$ is a function with real values in $[0,1] \subset \R$. A useful introduction to random graphs is \cite{frieze2015introduction}.

\begin{Conjecture}[BPR-conjecture]
\label{conj-bpr}
\ \\
For $\mathbb{G}(n, p(n))$ with $p(n)=1/2$,
the chromatic polynomial is almost complete.
\end{Conjecture}

They conjectured this 25 years ago, but today this seems unlikely.
As of today no almost complete graph polynomial is known.
In \cite[Chapter 14.3]{bk:BrouwerHaemers2012}  it is suggested that the characteristic
polynomial $char(G;x)$ might be almost complete.
Candidates for being almost complete multivariate polynomials are the Tutte polynomial, and if this is not the case, the
Edge Elimination polynomial $\xi(G;x,y,z)$ of \cite{averbouch2010extension} or the Bichromatic Polynomial from \cite{dohmen2003new}.

In contrast to Conjecture \ref{conj-bpr}
i{\em we} conjecture

\begin{Conjecture}
\label{conj-2}
\ \\
There is no Harary polynomial which is almost complete for $\mathbb{G}(n, 1/2)$.
\end{Conjecture}

\subsection{Weakly distinguishing graph polynomials}
Let $\mathbb{G}(n, p(n))$ be the random graphs of order $n$ where each edge has probability $p(n)$ for $0 \leq p(n) \leq 1$.
A good reference for random graphs is \cite{frieze2015introduction}.

In \cite{noy2003graphs} D. Welsh is credited for asking whether there are {\em naturally defined} graph polynomials
$F(G;x)$ such that almost all graphs are not $F$-unique. 
\begin{definition}
We say that a graph polynomial $F$ is {\em weakly distinguishing} if almost all graphs are not $F$-unique.
In other words 
$$
\lim_{n\rightarrow \infty} \dfrac{|U_F^{\cD}(n)|}{\cG(n)|}=0,
$$
or, equivalently, almost all graphs in $\mathbb{G}(n, 1/2)$ have an $F$-mate.
\end{definition}

M. Noy showed\footnote{
In \cite{noy2003graphs} the independence polynomial is called the stable set polynomial and denoted by $S(G;x)$.
} 
the following:
\begin{theorem}[M. Noy, 2003]
\label{pr:noy}
The independence polynomial $\rIND(G;x)$ is weakly distinguishing for $\mathbb{G}(n, 1/2)$.
\end{theorem}
Similar and more general results were also published in \cite{rakita2019weakly}.

Theorem \ref{pr:noy} does not address the question of the existence of $\rIND(G;x)$-unique graphs.
However, there are infinitely many $\rIND(G;x)$-unique graphs.

In \cite[Theorems 2.5-6]{beaton2019independence} we find proofs of the following:
\begin{theorem}
\label{th:ind-unique}
\ 
\begin{enumerate}[(i)]
\item
The graphs $K_n$ and $P_{2n+1}$ are $\rIND(G;x)$-unique.
\item
For every $k$
there is an even path $P_{2n(k)}$ with at least $k$ $\rIND(G;x)$-mates.
\end{enumerate}
\end{theorem}
There are several other families of $\rIND(G;x)$-unique graphs, such as the edge-less graphs,
or the graphs with only one edge.  Proofs can be given by direct calculations.

\begin{problem}
Find more families of $\rIND(G;x)$-unique 
graphs.
\end{problem}

More weakly distinguishing graph polynomials for $\mathbb{G}(n, 1/2)$, which are not Harary polynomials, were studied in \cite{rakita2019weakly}.
The analogue for hypergraph polynomials was discussed in \cite{makowsky2019p}.

\begin{problem}
\ 
\begin{enumerate}[(i)]
\item
Are there weakly distinguishing Harary polynomials for $\mathbb{G}(n, 1/2)$?
\item
Characterize $p(n)$ such that there are weakly distinguishing Harary polynomials for $\mathbb{G}(n, p(n))$
\end{enumerate}
\end{problem}

We shall discuss weakly distinguishing Harary polynomials for $\mathbb{G}(n, p(n))$ in Section \ref{se:wd} for $\cD = \cG$.
Theorem \ref{th:delta} in Section \ref{se:wd} shows the existence of weakly distinguishing Harary polynomials for $\mathbb{G}(n, d/n)$
with $d$ constant.

\subsection{Distinguishability index}
For a graph polynomial $F$, being almost complete and and weakly distinguishing are two extreme situations. 
We also want to look at intermediate notions.
Instead of looking at all graphs in $\cG$ we could restrict ourselves to a property $\cD$ and ask whether $F$ is
(almost) complete or weakly distinguishing on $\cD$. Instead of almost complete we can also consider the case where a certain proportion of
graphs  is $F$-unique on $\cD$.

More precisely, we look at the following:
Let $F$ be a graph polynomial, $\cD$ be a graph property.
$\cD(n)$ is the family of graphs in $\cD$ 
of order $n$ with vertices labeled $1,...n$. 
$U_F^{\cD}(n)$ is the family of $F$-unique graphs in $\cD$ of order $n$.
\begin{definition}
\ 
\begin{enumerate}[(i)]
\item
A graph polynomial $F$ is {\em $\alpha$-distinguishing on $\cD$} if 
$$
\lim_{n\rightarrow \infty} \dfrac{|U_F^{\cD}(n)|}{\cD(n)|}=\alpha.
$$
\item
A graph polynomial $F$ is {\em weakly distinguishing on $\cD$} if  $\alpha = 0$.
\item
$F$ is {\em almost complete on $\cD$} if $\alpha =1$.
The number $\alpha$ measures the distinguishing power of $F$ on $\cD$.
\\
We call $\alpha$ the {\em distinguishability index of $F$ on $\cD$}.
\end{enumerate}
\end{definition}

For $\mathbb{G}(n, p(n))$ the Distinguishability index can be defined similarly, but we do not need this in the sequel of the paper.

\begin{question}
Let $\chi_{\cP}$ be a Harary polynomial and let $\cD$ be a graph property.
Can we say anything interesting about the distinguishability index of $\chi_{\cP}$ on $\cD$?
\end{question}

\subsection{Comparability of Harary polynomials}

There are many graph polynomials which are d.p.-equivalent to the chromatic polynomial.
Similarly, for every graph property $\cP$ there are also many graph polynomials $F$
which are d.p.-equivalent to $\chi_{\cP}$.

\begin{proposition}
Let $\chi_{\cP} = \sum_i c^{\cP}_i(G) x_{(i)}$ and $F = \sum_i d^{\cP}_i(G) x_{(i)}$
\\
with $d^{\cP}_i(G) = f(i) c^{\cP}_i(G)$ where $f: \N^+ \rightarrow \N^+$ is any function.
\\
Then $\chi_{\cP}$ and $F$ are d.p.-equivalent.
\end{proposition}
\begin{proof}
It suffices to note that for any two graphs $G,H$ we have that 
\\
$\chi_{\cP}(G;x) = \chi_{\cP}(H;x)$ iff $F(G;x) = F(H;x)$.
\end{proof}

However, in the above proof, $F$ is not a Harary polynomial, if $f$ is not the identity.
Given two non-trivial graph properties $\cP, \cQ$, 
it is not clear under what conditions about $\cP$  and $\cQ$,
the two Harary polynomials $\chi_{\cP}$ and $\chi_{\cQ}$ are d.p.-comparable.

\begin{Conjecture}
\label{conj-3}
Given two non-trivial graph properties $\cP, \cQ$, 
then either $\chi_{\cP}$ and $\chi_{\cQ}$ are d.p-equivalent or d.p.-incomparable.
\end{Conjecture}
\ifnewpage
\newpage
\else \fi 

\section{The existence of $\chi_{\cP}$-mates}
\label{se:complete}
\ifmargin
\marginpar{\red File: 
\\
CM-complete.tex
\\
Last updated:
\\
16.12.25}
\else \fi 
In this section we present partial results supporting the Conjecture \ref{conj-1}.
We denote by $\cP(n)$ and  $\cP_<(n)$ the class of graphs in $P$ of order $n$, respectively of order 
at most $n$.
\\
We will frequently use the following lemma.
\begin{lemma}
\label{le:mates}
Let $m \in \N^+$ and $\cP$ and $\cQ$ be a two graph properties such that $\cP_<(m) = \cQ_<(m)$.
If  $G, H$  are two graphs of order at most $m$ then $G$ and $H$ are
$\chi_{\cP}$-mates iff they are $\chi_{\cQ}$-mates.
\end{lemma}
\begin{proof}
This follows immediately from Proposition \ref{pr:basic}.
\end{proof}

\subsection{Harary polynomials  of trivial graph properties}

A graph property is {\em trivial} if it is empty, finite or contains all finite graphs.
\begin{proposition}
Let $\cP$ be a trivial graph property and $\chi_{\cP}$ its Harary polynomial.
\begin{enumerate}[(i)]
\item
If $\cP = \emptyset$ the Harary polynomial  $\chi_{\emptyset}(G; k) =0$ for all $G$ and $k$.
\item
If $\cP = \cG$ consists of all graphs, its Harary polynomial of $G$ of order $n$ satisfies
$$\chi_{\cG}(G; k) = \sum_{i=0}^{n} S(n,i) k_{(i)}$$
where $S(n,k)$ are the Stirling numbers of the second kind. 
\item
Let $\cP$ be finite and $\cP(n)$ ($\cP_<(n)$) be consist of the graphs of $\cP$ of order \\
(at most) $n$.
\begin{enumerate}[(a)]
\item
If $\cP_<(2) = \emptyset$ the graphs $K_2$ and $E_2$ are mates.
\item
If $\cP_<(3) = \cG_<(3)$ the graphs $K_2$ and $E_2$ are mates.
\black
\end{enumerate}
Hence, Harary polynomials of trivial graph properties always have a pair of mates.
\end{enumerate}
\end{proposition}

\subsection{Almost all Harary polynomials have at least one pair of mates}
Recall that
$h_i^{\cP}(G)$ is the number ${\cP}$-colorings of $G$ with {\em exactly $i$ colors}.
$$
\chi_{\cP}(G;k) = \sum_i^n h_i^{\cP}(G) k_{(i)}.
$$
Obviously, in order to compute $\chi_{\cP}(G; k)$
it suffices to compute all the $h_i^{\cP}(G)$.

We first show that for many graph properties $\cP$ the weak form of Conjecture \ref{conj-1} is true.

\begin{proposition}
\label{pr:almost}
Let $\cP$ be a graph property and $G, H$ be two graphs of order at most $m$ which are $\chi_{\cP}$-mates.
\begin{enumerate}[(i)]
\item
Let $\cQ$ be a graph property such that $\cP(n) = \cQ(n)$ for all $n \leq m$.
Then $G$ and $H$ are also $\chi_{\cQ}$-mates. 
\item
For almost all Harary polynomials $G_{\cP}$ and $H_{\cP}$ also mates.
\end{enumerate}
\end{proposition}
\begin{proof}
(i):
We use that $G_{\cP}, H_{\cP}$ are of order at most $m$ and $\cP(n) = \cQ(n)$ for all $n \leq m$.
\\
(ii):
There are only finitely many sets of graphs (up to isomorphisms) of order at most $m$.
But there are continuum many non-trivial graph properties.
\end{proof}


\begin{proposition}
Let $n_0 \in \N^+, n_0 \geq 2$.
\\
Assume $P_{\leq}(n_0)$ does {\bf not} contain any edge-less graphs. 
For all $n \leq n_0, k \in \N^+$ we have  $\chi_P(E_n, k) = 0$. 
Hence $E_j$ and $E_{j'}$ are mates for all $j \neq j' \leq n_0$.
\end{proposition}
\begin{proof}
We compute $h_i^P(E_{n})$ for $n \leq n_0, i \in [n]$:
$$
h_i^P(E_n) = 0,   n \leq n_0, i \in [n]
$$
\end{proof}

\subsection{Harary polynomials of Compton-Gessel classes}

A graph property $\cP$ is a {\em Compton-Gessel class} if it is closed under taking disjoint unions and connected components.
Compton-Gessel classes were introduced and studied first in \cite{Compton-PhD,compton1983some,gessel1984combinatorial}.

\begin{proposition}
\label{pr:multiplicative}
Let $\cP$ be a graph property.
The Harary polynomial $\chi_{\cP}$ is multiplicative iff $\cP$ is a Compton-Gessel class.
\end{proposition}
\begin{proof}
Assume $\cP$ is a Compton-Gessel class. 
Let $A \subseteq V(H_1) \sqcup V(H_2)$ be a color class of $f$ and $A_i = A \cap V(H_i), i = 1,2$.
If $f$ is a $\cP$-coloring of  $H =(H_1 \sqcup H_2)$ with at most $k$ colors, then $H(A) \in \cP$ and each $A_i$ 
induces a graph which is a disjoint union of connected components of $H(A)$.
Therefore $H(A_i) \in \cP$. Hence the restriction of $f$,  $f_i = f|_{A_i}a,$ is a $\cP$-coloring of each $H(A_i)$ with at most $k$ colors.

On the other hand, if each  $f_i$ is a $\cP$-coloring of $H_i$, the disjoint union $f_1 \sqcup f_2$ is a 
$\cP$-coloring of $H_1 \sqcup H_2$. Therefore
$$
\chi_{\cP}(H_1 \sqcup H_2, k) = \chi_{\cP}(H_1, k) \cdot \chi_{\cP}(H_2, k).
$$

Now assume  $\chi_{\cP}$ is multiplicative. So, every $\cP$-coloring  $f$ of $(H_1 \sqcup H_2)$
is of the form $f =f_1 \sqcup f_2$. Therefore $\cP$ is closed under disjoint unions and connected components, hence a Compton-Gessel class.  
\end{proof}

The following are Compton-Gessel classes:
\begin{examples}
\begin{enumerate}[(i)]
\item
The class $\cE$ of edge-less graphs.
\item
The planar graphs.
\item
Let $H$ be a connected graph. The class $Forb_{ind}(H)$ of graphs with no induced subgraph isomorphic to $H$ is a Compton-Gessel class.
It is also hereditary.
\item
The class of graphs which are disjoint unions of triangles is a Compton-Gessel class but it is not hereditary.
\item
Let $\cP_0$ be a finite hereditary class and $\cE$ the edge-less graphs. $\cP_0 \cup \cE$ is non-trivial hereditary but not Compton-Gessel.
\item
$TW(k)$, the graphs of tree-width at most $k$,  and also
$CW(k)$, the graphs of clique-width at most $k$, are Compton-Gessel classes.
\item
$\cP_t$, where $\cP_t$ is the graph property such that $H \in \cP_t$ iff
the connected components of $H \in \cP_t$ have order at most $t$.
\item
The intersection of two Compton-Gessel classes is again a Compton-Gessel class.
Every graph property $\cP$ is a subset of a minimal Compton-Gessel class.
\end{enumerate}
\end{examples}

The following follows from the definitions: 

\begin{proposition}
Let $\cP$ a Compton-Gessel class of graphs and let $G_1, G_2$ be $\chi_{\cP}$-mates.
For every graph $G \in \cP$ the graphs $G_1 \sqcup G$ and $G_2 \sqcup G$ are also $\chi_{\cP}$-mates.
\end{proposition}

\subsection{$\chi_{\cP}$-mates in $\cP$ vs $\chi_{\cP}$-mates not in $\cP$}

We shall see in the sequel that it is easier to formulate general theorems about the existence of $\chi_{\cP}$-mates in $\cP$.
The proofs that certain graphs are chromatically equivalent and are not in $\cE$ does not seem to generalize in an obvious way.
They all depend ultimately on more or less sophisticated uses of the deletion/contraction relations of the chromatic polynomial,
see \cite[Chapter 3]{bk:DongKohTeo2005}.

A trivial example of $\chi_{\cP}$-mates, which are not in $\cP$, can be constructed as follows.

\begin{proposition}
Let $\cP$ be such that its complement is hereditary and let \\
$G \not \in \cP$.
Then $\chi_{\cP}(G; k) = 0$ for all $k$.
\\
In particular, any two non-isomorphic graphs $G_1, G_2 \not \in \cP$ are $\chi_{\cP}$-mates.
\end{proposition}
The proofs are left to the reader.

\ifnewpage
\newpage
\else \fi 

\section{Using heredity}
\label{se:hereditary}
\ifmargin
\marginpar{\red File: 
\\
CM-hereditary.tex
\\
Last updated:\\
16.12.25}
\else\fi 

\subsection{Hereditary graph properties}

Let $\cP$ be a graph property of loop-free simple graphs.
$\cP$ is {\em hereditary} if it is closed under induced subgraphs.
It is monotone if it is closed under subgraphs.
$\cP$ is {\em minor closed} if it is is closed under graph minors.
For two graphs $G', G$ we write 
$G' \subset_s G$
($G' \subset_i G$, $G' \subset_m G$) if $G'$ is a subgraph (induced subgraph, minor) of $G$.
$\cP$ is {\em partially monotone (hereditary, minor closed)} if there is a property 
$\cP_0 \subseteq \cP$ such that for all $G \in \cP_0$
with $G' \subset_s G$ ($G' \subset_i G$,  $G' \subset_m G$) we have $G' \in \cP_0$.

\begin{facts}[Hereditary properties]
\ 
\begin{enumerate}[(i)]
\item
If $\cP$ is minor closed it is monotone.
\\
If $\cP$ is monotone it is also hereditary.
\item
Every non-empty hereditary property contains $E_1$.
\end{enumerate}
\end{facts}

We denote by $\cP^*(n)$ the number of isomorphism classes of graphs of order $n$. Note that $\cP^*(n) < \cP(n)$.
We count in $\cP(n)$ the labeled graphs, and in $\cP^*(n)$ the unlabeled graphs.
We define the {\em labeled speed} of $\cP$, by $s_{\cP}(n) =|\cP(n)|$ and
the {\em unlabeled speed} of $\cP$, by $s_{\cP}^*(n) =|\cP^*(n)|$.

\begin{theorem}[\cite{balogh2009unlabelled}]
\label{th:BBZS-speed}
Let $\cP$ be hereditary. Then $s_{\cP^*}(n)$ satisfies one of the following:
\begin{enumerate}[(i)]
\item
$s_{\cP}(n)$ is ultimately constant with $s_{\cP}(n) \in \{0, 1, 2\}$.
\\
$s_{\cP}(n) =0$ for large enough $n$ iff $\cP$ is finite.
\item
There are constants $k \in \N^+$ and $c \in \Q$ such that for all $n$ 
$$
s_{\cP}(n) = cn^k + O(n^{k-1}).
$$
\item
Ultimately $s_{\cP}(n) \geq S(n)$ where $S(n)$ is the number of partitions of a set with
$n$ indistinguishable elements into nonempty subsets.
\end{enumerate}
\end{theorem}

\begin{corollary}
\label{co:hereditary}
If $s_{\cP}(n)$ is not ultimately constant $s_{\cP}(n) = 0, 1$,
then $s_{\cP}(n) \geq 2$ for sufficiently large $n$.
In other words, for sufficiently large $n$ there are at least two non-isomorphic graphs in  $\cP(n)$. 
\end{corollary}

Recall that $\cG$ is the class of all finite graphs up to isomorphism.
\begin{proposition}
\label{pr:hereditary-small}
Let $\cP$ be hereditary.
Then 
$$
\lim_{n\rightarrow \infty} \dfrac{|\cP(n)|}{|\mathcal{G}(n)|} =0.
$$
\end{proposition}
\begin{proof}
As $\cP$ be hereditary there exists a set of graphs $\cH$ such that $\cP$ is the set of graphs which do not have a graph $H \in \cH$ as
an induced subgraph. Let $Forb(\cH)$ the set of these graphs and $Forb(H)$ the same for a single $H$.  
For any graph $H \in \cH$ it is well known, \cite{frieze2015introduction},  that
$$
\lim_{n\rightarrow \infty} \dfrac{|Forb(H)(n)|}{|\mathcal{G}(n)|} =0.
$$
hence
$$
\lim_{n\rightarrow \infty} \dfrac{|Forb(\cH)(n)|}{|\mathcal{G}(n)|} =0.
$$
\end{proof}
Early results for the labeled case
can be found in \cite{alekseev1992range,balogh2000speed}.
A good survey of the rates of growth of graph properties  (labeled and unlabeled) is \cite[Sections 2.2-3]{klazar2010some}.

\subsection{Finding many  $\chi_\cP$-mates}
In this section we discuss how to use (partial) heredity in order to find $\chi_\cP$-mates in $\cP$.

\begin{proposition}
\begin{enumerate}[(i)]
Let $\cP$ be a non-trivial hereditary graph property. 
\item
If $\cP$ contains a graph $G$ of order $n$ which is neither a clique nor the complement of a clique
which has an edge $e$, and two vertices which are not an edge,  then $K_2$ and $E_2$ are also in $\cP$ and they are $\chi_{\cP}$-mates.
\item
If for some $n$ there are at least two non-isomorphic graphs $G, G' \in \cP(n)$, the graphs $G,G'$ are $\chi_{\cP}$-mates.
\item
The same holds if $\cP$ is monotone or minor-closed.
\end{enumerate}
\end{proposition}
\begin{proof}
(i): By Lemma \ref{le:mates} 
we have to check $\chi_{\cP}(K_2,k)$ and $\chi_{\cP}(E_2,k)$ for $k=1,2$.
\\
$\chi_{\cP}(K_2,1) = \chi_{\cP}(E_2,1) =1$ by Lemma \ref{le:1}.
\\
$\chi_{\cP}(K_2,2) = \chi_{\cP}(E_2,2) =1$ since $E_1 \in \cP$.

(ii): We have to check $\chi_{\cP}(G;k)$ and $\chi_{\cP}(G';k)$ for $k \in [n]$.
\\
If $G \in \cP$ and $\cP$ is hereditary, every subset $X$ of $V(G)$ induces a graph $G[X] \in \cP$.
We conclude that
$$
\chi_P(G;x) = \chi_P(G';x) = \sum_i^n S(i,k) x_{(i)}.
$$
(iii): Follows from (ii).
\end{proof}

Using the above proposition we get immediately:

\begin{theorem}
\label{th:hereditary}
Let $\cP$ be a non-trivial hereditary graph property such that $s_{\cP}(n)$ is not ultimately constant $s_{\cP}(n) = 0, 1$.
Then there are infinitely many pairs of graphs which are $\chi_{\cP}$-mates.
\end{theorem}

\begin{examples}
If $\cP$ is one of the following graph properties then
$\chi_{\cP}$ has infinitely many pairs of graphs which are $\chi_{\cP}$-mates.
\begin{enumerate}[(i)]
\item
If $\cP = \cE$, i.e., $\cP$ consists of all edge-less graphs, or $\cP$ consists of all complete graphs,
$\cP$ is hereditary but $s_{\cP}(n) = 1$ for all $n \in \N^+$.
This gives the chromatic polynomial or the related  adjoint polynomial from the Examples \ref{ex:harary}.
We have seen that $\chi$ has infinitely many pairs of graphs which are $\chi_{\cP}$-mates.
However, they are not in  $\cE$. Actually no graph in $\cE$ has a $\chi$-mate.
\item
The $mcc_t$-colorings from Example \ref{ex:harary}. 
Here $\cP = \cP_t$ is the graph property such that $H \in \cP$ if all the connected components of $H$ 
have order at most $t$.  $\cP_t$ is hereditary and $s_{\cP}(n) \geq  2$ for $t \geq 2$.
Any two graphs in $\cP_t$ of the same order are $\chi_{mcc_t}$-mates.
\end{enumerate}
\end{examples}


\subsection{Finding $\chi_P$-mates for $\cP$ a partially hereditary property}

The proof of Theorem \ref{th:hereditary} actually shows the following:

\begin{proposition}
Let $\cP$ be a graph property 
and  $\cP_0 \subseteq \cP$ be hereditary and possibly finite.
Any two non-isomorphic  graphs in $\cP_0$ of the same order are $\chi_{\cP}$-mates.
\end{proposition}
Together with Proposition \ref{pr:almost} we get:
\begin{corollary}
Fix $k \in \N^+$.
Almost all Harary polynomials have at least $k$ mates.
\end{corollary}

\ifnewpage
\newpage
\else \fi 

\section{Finding mates for graph polynomials as generating functions}
\label{se:generating}
\ifmargin
\marginpar{\red File:
\\
CM-newmates.tex
\\
Last updated:
\\
16.12.25}
\else \fi 
In this section we apply the techniques of Section \ref{se:complete}
to graph polynomials as generating functions.

For a graph property $\cP$ let $c_{\cP}(G, i)$ be the number of induced subgraphs of $G \in \cP$ of order $i$.
We look at graph polynomials of the form
$$
F_{\cP}(G; x) = \sum_i  c_{\cP}(G; i) x^i.
$$
They are generating functions for the number of induced subgraphs in $\cP$.

\begin{theorem}
Let $\cP$ be hereditary
and $G, H \in \cP$ both of order $n$.
\\
Then $F_{\cP}(G; x) = F_{\cP}(H; x)$.
Hence, if $G \ncong  H$, they are $F_{\cP}$-mates.
\end{theorem}

\begin{proof}
Assume $V(G)=V(H)= [n]$. Let $A \subseteq [n]$. Both $G(A)$ and $H(A)$ are in $\cP$.
Hence, $c_{\cP}(G, i) = c_{\cP}(H, i)$ for all $i \in [n]$.
\end{proof}

We are left with the case, where $\cP$ is non-trivial and hereditary, \\
and $|\cP(n)| = 1$ for all $n$.

\begin{proposition}
Assume now that $\cP$ is non-trivial and hereditary, \\
and $|\cP(n)| = 1$ for all $n$. 
\begin{enumerate}[(i)]
\item
$\cP = \cE$, i.e., it consists only of edge-less graphs or only of complete graphs.
\item
$F_{\cP}(G; x)$ is the independence polynomial $\rIND(G;x)$ or the clique polynomial $\mathrm{Cl}(G;x)$.
The clique polynomial is the generating function $\mathrm{Cl}(G;x) = \sum_i c_i(G) x^i$ where $c_i(G)$ counts the number of induced clique of order $i$.
\\
Note that $\rIND(G;x)= \mathrm{Cl}(\bar{G};x)$ where $\bar{G}$ be the complement graph of $G$.
\item
No $G, H \in \cP$ can be $F_{\cP}$-mates. 
\end{enumerate}
\end{proposition}
\begin{proof}
(i): If $G \in \cP$ has an edge $e =(u,v)$, $K_2 = (\{u,v\}, \{e\})$ and \\
$E_2 = (\{u,v\}, \emptyset)$ are in $\cP$,
hence $|\cP(2)| = 2$, a contradiction.
\\
(ii) and (iii): By definition.
\\
(iv): 
For $G \in \cP$ we have $c_{\cP}(G, i) =1$ iff $G$ is of order $i$.
If $G$ and $H$ have the same order, they are isomorphic, and otherwise their polynomials differ.
\end{proof}

Mates for the independence polynomial have been studied in the literature,
\cite{chism2009independence,beaton2019independence}.
Let $D_n$ be the graph of order $n$  which consists of a triangle with a tail (path) of 
$n-3$ vertices attached to a vertex of the triangle. 
We have seen in Theorem  \ref{pr:noy} that almost all graphs have an $\rIND(G;x)$-mate.
In the literature explicit $\rIND(G;x)$-mates are rarely given. 

\begin{theorem}[\cite{chism2009independence,beaton2019independence}]
$C_n$ and $D_n$ are  $\rIND(G;x)$-mates.
and $\bar{C_n}$ and $\bar{D_n}$ are $Cl(G;x)$-mates.
\end{theorem}

Let $A_{\cP} = \{ n \in N^+: |\cP(n)| \geq 2 \}$.
For graph properties $\cP$ of  loop-free simple graphs $1 \not\in A_{\cP}$.
We summarize the above.

\begin{theorem}
\label{th:gen-hereditary}
Let $\cP$ be non-trivial and hereditary. 
\begin{enumerate}[(i)]
\item
If $A_{\cP} \neq \emptyset$ is finite with $n$ elements there are at least $n$ pairs of graphs $G,H$ which are $F_{\cP}$-mates.
\item
If $A_{\cP}$ is infinite, there are infinitely many pairs of graphs $G,H$ which are $F_{\cP}$-mates.
\item
If $A_{\cP} = \emptyset$, there are infinitely many pairs of graphs $G,H$ which are \\
$F_{\cP}$-mates.
\end{enumerate}
\end{theorem}

\ifnewpage
\newpage
\else \fi 

\section{Comparing Harary polynomials}
\label{se:compare}
\ifmargin
\marginpar{\red File: \\
CM-compare.tex
\\
Last revised:
\\
16.12.25}
\else \fi 

In this section we 
we study d.p.-incomparability of Harary polynomials.

\subsection{How to prove d.p.-incomparability?}

Let $F_1, F_2$ be two graph polynomials (graph invariants).
The following proposition follows from the definitions from Section \ref{se:intro}.
\begin{proposition}
Two graph parameters
$F_1$ and $F_2$ are d.p.-incomparable iff there are graphs $H_1, H_2, H_1', H_2'$ such that
$$
F_1(H_1) = F_1(H_1') \text{  and  } F_2(H_1) \neq F_2(H_1')
$$
and
$$
F_1(H_2) \neq F_1(H_2') \text{  and  } F_2(H_2) = F_2(H_2')
$$
In other words $H_1, H_1'$ are $F_1$-mates but not $F_2$-mates,
and 
$H_2, H_2'$ are $F_2$-mates but not $F_1$-mates.
\end{proposition}

In \cite{makowsky2019logician} it is shown that the chromatic polynomial $\chi(G;x)$, the characteristic polynomials $\mathrm{char}(G;x)$ and $\mathrm{Lap}(G;x)$, 
and the independence polynomial $\rIND(G;x)$ are pairwise d.p.-incomparable.

\subsection{Harary polynomials d.p.-incomparable to $\chi$}
\begin{theorem}
\label{th:dp-1}
\ \\
There are uncountably many Harary polynomials which are d.p.-incomparable to the chromatic polynomial.
\end{theorem}
\begin{proof}
Let $G  \not \sim H$ be two $\chi$-mates, 
hence $\chi(G;x) = \chi(H; x)$.
\\
Let $\cP_H = \{H\}$ be the trivial graph property consisting of $H$ only.
\\
So $\chi_{\cP_H}(G;x) = 0 \neq \chi_{\cP_H}(H; 1) = 1$ and $G$ and $H$ are not $\chi_{\cP_H}$-mates.
\\
As the chromatic polynomial is multiplicative, we get
$$\chi(G \sqcup G; x) = \chi(G; x) \cdot \chi(G; x).$$
\\
If $\chi(G; x) \neq 1$ for some value of $x$ we have $\chi(G \sqcup G; x) \neq \chi(G; x)$ and 
$G$ and $G \sqcup G$ are not $\chi$-mates.
\\
On the other hand we have
$\chi_{\cP_H}(G;x) = \chi_{\cP_H}(G \sqcup G; x) = 0$.
Hence $G$ and $G \sqcup G$ are $\chi_{\cP_H}$-mates.
\item
As noted in Example
\ref{chi-mates}, all trees of the same order are $\chi$-mates.
Hence, there are infinitely many pairs $G, H$ which are $\chi$-mates.

As noted in Lemma \ref{le:mates},
computing $\chi$ and $\chi_{\cP_H}$ for $G,H$ and $G \sqcup G$
does not depend on the graphs in $\cP_H$ which are larger than $G,H$ and $G \sqcup G$.
Therefore, there are uncountably many graph properties $\cP$  for which the same argument works.
Hence, there are uncountably many Harary polynomials which are d.p.-incomparable to the chromatic polynomial.
\end{proof}

\subsection{Mutual incomparability of Harary polynomials}

Let $\{G_i: i \in \N\}$ be an enumeration of all finite graphs up to isomorphism.

Let $P_i = \{G_i\}$ be the trivial graph property consisting of the single graph $G_i$ and
$\chi_i(G; x)$ be the Harary polynomial defined by $P_i$.

\begin{proposition}
\label{pr:singleton}
Assume both $G_i$ and $G_j$ are connected and of the same order.
Then the Harary polynomials $\chi_i(G; x)$ and $\chi_j(G; x)$ are d.p.-incomparable.
\end{proposition}
\begin{proof}
We look at the four graphs $G_i, G_j, G_i \sqcup G_i, G_j \sqcup G_j$ and compute $\chi_i$ for each of them.
Analogously, we can also compute $\chi_j$.
\begin{description}
\item[$\chi_i(G_i; k)$]
$\chi_i(G_i; 1) = 1$ and $\chi_i(G_i; k) = 0$ for $k \geq 2$ because each block of a partition into more than one 
part is too small to induce $G_j$.
\item[$\chi_i(G_j; k)$]
$\chi_i(G_j; 1) = 0$
and $\chi_i(G_i; k) = 0$ for $k \geq 2$ because each block of a partition into more than one 
part is too small to induce $G_i$.
\item[$\chi_i(G_j \sqcup G_j; k)$]
$\chi_i(G_j \sqcup G_j; k) = 0$ for all $k \geq 1$.
\\
For $k =1$  we use $\chi_i(G; 1) = 1$ iff $G$ is isomorphic to $G_i$.
\\
For $k \geq 2$  we use that both $G_i$ and $G_j$ are connected and of the same order.
\item[$\chi_i(G_i \sqcup G_i; k)$]
$\chi_i(G_i \sqcup G_i; k) \neq 0$ for $k=2$.
\end{description}
We conclude that $G_j$ and $G_j \sqcup G_j$ are $\chi_i$-mates but not $\chi_j$-mates.
Analogously $G_i$ and $G_i \sqcup G_i$ are $\chi_j$-mates but not $\chi_i$-mates.
Hence $\chi_i$ and $\chi_j$ are d.p.-incomparable.
\end{proof}

The graph properties $\cP_i$ underlying $\chi_i(G; k)$ are singletons.
Using Lemma \ref{le:mates}
we can modify the proof of Proposition \ref{pr:singleton} above such that these properties are infinite
and even disjoint.

\begin{theorem}
\label{th:dp-2}
There are uncountably many disjoint non-trivial pairs of graph properties $\cP, \cQ$ such that
$\chi_{\cP}$ and $\chi_{\cQ}$ are d.p.-incomparable.
\end{theorem}


Before we prove or disprove Conjecture \ref{conj-3}, we might try to answer the following questions:
\begin{question}
\begin{enumerate}[(i)]
\item
Let
$\cP \subseteq \cQ$ be two graph properties. What can we say in general about 
$\chi_{\cP}  <_{d.p} \chi_{\cQ}$?
\item
Let $\cP$ and $\cQ$ be two distinct graph properties.
Under which conditions are
$\chi_{\cP}$  are $\chi_{\cQ}$ d.p.-incomparable?
\item
Let $\cP$ and $\cQ$ be any two disjoint graph properties.
Is it always the case that $\chi_{\cP}$  are $\chi_{\cQ}$ d.p.-incomparable?
\end{enumerate}
\end{question}

\ifnewpage
\newpage
\else \fi 

\section{Weakly distinguishing Harary polynomials}
\label{se:wd}
\ifmargin
\marginpar{\red File: 
\\
CM-wd.tex
\\
Last updated:
\\
24.12.25}
\else \fi 

Recall that a graph polynomial $F$ is trivial if any two pairs of graphs are $F$-mates.
$F$ is weakly distinguishing if almost all graphs are not $F$-unique.

\subsection{Harary polynomials $\chi_{\cP}$ with very few $\chi_{\cP}$-unique graphs}

The Harary polynomial $\chi_{\emptyset}$ is trivial. $\chi_{\emptyset}(G;k) = 0$ for all graphs $G$ and $k \in \N^+$.
The Harary polynomial $\chi_{\cG}$ with $\cG$ all finite graphs satisfies $\chi_{\cG}(H;k) = S(n,k)$ for every $H$ of order $n$,
where $S(n,k)$ is the Stirling number of the second kind. 
We conclude that $K_1$ is the only  $\chi_{\cG}$-unique graphs and any two graphs of the same order are $\chi_{\cG}$-mates.

\begin{proposition}
There are uncountably many graph properties $\cP_A$ such that
$K_1$ is the only  $\chi_{\cG}$-unique graph and any two graphs of the same order are $\chi_{\cG}$-mates.
\end{proposition}
\begin{proof}
Let $A \subseteq \N^+$ and let $\cP_A$ be the graph property which consists of all graphs whose order has a size in $A$.
There are uncountable many such sets $A$.

Let $S_A(n,k)$ be the number of partitions of $[n]$ into $k$ blocks, each of which has a size in $A$.
They are the {\em restricted Stirling numbers} of the second kind first studied in \cite{broder1984r}.
We get $\chi_{\cP_A}(H;k) = S_A(n,k)$ for every $H$ of order $n$.
\end{proof}

All the examples so far are obviously weakly distinguishing for $\mathbb{G}(n, 1/2)$. 

\subsection{The Harary polynomials $F_r(G;k)$}
Let $\bigcirc_r$ the class of graphs which consists of disjoint unions of copies of $C_r$ and 
Let $mC_r$ be the graph which consists of disjoint unions of $m$ circles $C_r$.
$G \in \bigcirc_r$ iff  $G \cong mC_r$ for some $m \in \N^+$.

$F_r(G;k) = \chi_{\bigcirc_4}(G;k)$
be the partition Harary polynomial where in each partition of $V(G) = V_1, \ldots , V_k$ into $k$ blocks
each induced subgraph $G[V_i]$ is in $\bigcirc_r$.

Let $C(r,\ell)$ with $\ell \leq r$ be the class of graphs which do not contain exactly $\ell$  vertices adjacent to the same two cycles $C_r$.

\begin{proposition}
\label{pr:8-1}
\ 
\begin{enumerate}[(i)]
\item
$F_r$ is  is multiplicative for $r \geq 3$.
\item
Let $G_1, G_2$ be $F_r$-mates and $H$ be any graph. Then $G_1 \sqcup H$ and $G_2 \sqcup H$ are also $F_r$-mates.
\item
If $F_r(G;x)$ does not vanish, i.e., is not identically $0$, for $G$, every $v \in V(G)$ is part of a cycle $C_r$.
\item
If, furthermore, $G \in Forb(C(r,1)$, every $v \in V(G)$ is part of exactly one cycle $C_r$.
\item
If $G$ is a graph of order $n \neq 0 \pmod{r}$
the polynomial  $F_r(G;k)$ vanishes.
\end{enumerate}
\end{proposition}
\begin{proof}
\ \\
(i): $\bigcirc_r$ is 
is a Compton-Gessel class, i.e., it is 
closed under taking components and disjoint unions, 
so by Proposition \ref{pr:multiplicative}
$F_r(G;k)$ is multiplicative.
\\
(ii): Follows from (i).
\\
(iii): Let $F_r(G;k_0) \neq 0$. So $V(G)$ can be partitioned into $k_0$ blocks, where each block consists of disjoint unions of cycles $C_r$.
\\
(iv): Assume  $F_r(G;x)$ is not identically $0$ for $G$ and $G \in Forb(C(r,1))$. No vertex of $G$ can be part of two induced cycles.
\\
(v): This follows from (iv).
\end{proof}

Let $D_{r,m}$ be the graph which is the disjoint union of $m$ copies of $C_r$.
Let $\bigcirc_r^+$ be the class of graphs 
of the form $(V = [rm], E_m \sqcup E)$ such that $([3m], E_m) \in \bigcirc_r$.
We denote these graphs by $D_{r,m}^E$ and $D_{r,m}$ is the graph $D_{r,m}^E$  with $E = \emptyset$.

\begin{proposition}
\label{pr:unique}
\ 
\begin{enumerate}[(i)]
\item
We have
$F_r(D_{r,m}^E, k) = \sum_i^m S(m,i) k^i$ iff $E = \emptyset$.
\item
Let $G$ be a graph with
$F_r(G, k) = \sum_i^m S(m,i) k^i = F_r(D_{r,m}, k)$.
\\
Then $G \cong D_m$, hence $D_m$ is $F_r$-unique. 
\item
All graphs in $\bigcirc_r$ are $F_r$-unique, hence for each $r \geq 3$ there are infinitely many $F_r$-unique graphs.
\end{enumerate}
\end{proposition}

\begin{proof}
\ 
(i): $D_{r,m}^E$ consists of $m$ copies of $C_r$ with additional edges in $E$.
If $E=\emptyset$ every partition of these $m$ cycles into  $i$ groups of disjoint $C_m$'s is a $\bigcirc_r$-partition.
If $E \neq\emptyset$ at least one of these partitions is not a $\bigcirc_r$-partition.
\\
(ii): $G$ and $D_{r,m}$ have the same order $rm$. Furthermore, 
$$
F_r(D_{r,m}; 1)=F_r(G; 1) = 1.
$$ 
But $F_r(G; 1) = 1$ iff 
$G \in \bigcirc_r$, hence $G \cong D_{r,m}$.
\\
(iii)
This follows from (ii).
\end{proof}

\begin{proposition}
\label{pr:vanish}
\end{proposition}

\begin{theorem}
\label{th:delta}
\ 
\begin{enumerate}[(i)]
\item
$F_r$ is weakly distinguishing  
in $\mathbb{G}(n, d/n)$ 
for $r \geq 3$. 
\item
$F_r$ is weakly distinguishing  
in $\mathbb{G}(n, 1/2)$ 
for $r \geq 4$. 
\end{enumerate}
\end{theorem}


Before we prove Theorem \ref{th:delta} in Subsections \ref{sse:proof1} and \ref{sse:proof2}
we collect some material about random graphs, see \cite{spencer2013strange,frieze2015introduction}.

\subsection{Digression on random graphs}
\ 
From the $0/1$-law of first order logic, \cite{spencer2013strange} we get immediately:
\begin{proposition}
\label{le:0/1}
Let $H$ be a fixed graph. Almost all graphs  in $\mathbb{G}(n, 1/2)$ contain $H$ as an induced subgraph.
Hence, in $\mathbb{G}(n, 1/2)$ almost no graphs are in $Forb(H)$. 
\end{proposition}

From \cite{frieze2015introduction} we collect the following:
\begin{proposition}[{\cite[Theorem 3.5]{frieze2015introduction}}]
\label{le:degree-p-constant}
In  $\mathbb{G}(n, p)$ with $p$ constant, the maximum degree is $O(n)$.
\end{proposition}

This contrasts with the following:
\begin{proposition}[{\cite[Theorem 3.4]{frieze2015introduction}}]
\label{le:degree}
Almost all graphs  in $\mathbb{G}(n, d/n)$ with $d$ constant have  maximal degree $D \approx  \frac{\log n}{\log\log n}$.
\end{proposition}


For $p(n)= d/n$ with $d$ constant however, the following can be used.
\begin{proposition}[{\cite[Exercise 1.4.5]{frieze2015introduction}}]
\label{le:d/n}
Let $G$ be in $\mathbb{G}(n, d/n)$ with $d$.
\begin{enumerate}[(i)]
\item
With high probability (w.h.p.)
no vertex of $G$ belongs to more than one triangle $C_3$.
\item
Let $r \geq 4$.
With high probability
no vertex belongs to more than one circle $C_r$.
\end{enumerate}
\end{proposition}
\begin{proof}
This is an easy first moment argument.
\\
(i): We have to show that w.h.p. $G$ does not contain  two triangles with one or two common vertices,
It is enough to show that w.h.p. high probability $G$ does not contain a $C(3,1)$-graphs, i.e. a graph which consists of two triangles with
exactly one common vertex.
\\
There are $O(n^5)$ potential $C(3,1)$-graphs, and each of them appears in $G(n, d/n)$ with probability ${(\frac{d}{n}})^6$.
Hence the expected number of $C(3,1)$-graphs is 
$O(n^5) \cdot (\frac{d}{n})^6 = O(\frac{d^6}{n}) = o(1)$.
Therefore by Markov’s inequality, almost all graphs are in $Forb(C(3,1))$.
\\
(ii): We have to show that w.h.p. $G$ does not contain  two cycles $C_r$ with any number of common vertices.
Again, it is enough to show that w.h.p. height probability $G$ does not contain a graph which consists of two cycles $C_r$ with
exactly one common vertex. We denote such a graph with $C_r \sqcup_1 C_r = C(r,1)$
\\
The proof is a modification of the proof of (i).
There are $O(n^{2r-1})$ potential $C(r,1)$-graphs, and each of them appears in $G(n, d/n)$ with probability $(d/n)^{2r}$.
Hence the expected number of $C(r,1)$-graphs is $O(n^{2r-1}) \cdot (\frac{d}{n})^{2r} = O(\frac{d^{2r}}{n}) = o(1)$.
\end{proof}

\begin{remark}
In \cite[Exercise 1.4.5]{frieze2015introduction} it is only stated for triangles $C_3$.
\end{remark}

\subsection{Proof of Theorem \ref{th:delta}(i)}
\label{sse:proof1}
\ 
\begin{fact}
\label{step:8-1}
$\bigcirc_r \subset Forb(C(r,1))$.
\end{fact}

Let $\cQ$ be the class of graphs which are either graphs of order $n \neq 0 \pmod{r}$ or are $G$ is $C_r$-free.
and $\Q^* = \cQ \cap Forb(C(r,1)))$.
\begin{fact}
\label{step:8-2}
Assume $G$ and $G'$ are both in $\cQ$.
Then $F_r(G; k) = 0$ is identically $0$ and $G$  and $G'$ are $F_r$-mates.
\end{fact}

Let $\bigcirc_r^+$ be the class of graphs 
of the form $(V = [rm], E_m \sqcup E)$ such that $([rm], E_m) in \bigcirc_r$.
Let $\bigcirc_r^*$ be $\bigcirc_r^+ \cap Forb(\bigcirc_r)$.

\begin{fact}
\label{2-edge-connected}
Let $G \in \bigcirc_r^*$ be at least $2$-edge-connected, and let $e$ and edge of $G$ which is not part of a cycle $C_r$.
then $F_r(G;k) =F_r(G_{-e};k)$ for all $k \in \N^+$, hence $G$ and $G_{-e}$ are $F_r$-mates.
\end{fact}
\begin{proof}
$F_r(G;k)$ depends only on which triangles of $G$ are connected.
\end{proof}

\begin{fact}
\label{bridge}
Let $G \in \bigcirc_r^*$ be with an edge $b=(u,v)$, $u \in V(G_1)$ and $v \in V(G_2)$  which is a bridge connecting subgraphs $G_1$ and $G_2$ both in $\bigcirc_r^*$.
Then there exists vertices $u' \in V(G_1)$ and $v' \in V(G_2)$ such that adding the edge $d=(u',v')$ does not generate a circle $C_r$.
\end{fact}
\begin{proof}
As $b$ is a bridge, there are no edges in $G$ between $G_1$ and $G_2$. Choose $u' \in V(G_1),  u' \neq u$ and $v' \in V(G_2), v' \neq v$.
Then adding $b' = (u',v')$ to $G$ cannot generate a circle $C_r$.
\end{proof}

\begin{proposition}
\label{pr:8-2}
Every graph $G \in \bigcirc_r^*$ has a $F_r$-mate.
\end{proposition}
\begin{proof}
Without loss of generality (Proposition \ref{pr:8-1}(ii)) we can assume hat $G$ is connected.
\\
Assume $G$ is at least $2$-edge-connected. There are $d,e$ two edges of $G$ such that neither $d$ nor $e$ are part of a triangle.
By Fact \ref{bridge}
$G_{-d}$ and $G_{-e}$ are $\chi_{\Delta}$-mates.
\\
If $G$ has a bridge $b=(u,v)$, we can apply Fact \ref{bridge} and choose $b'=(u',v')$. Let $G'$ be obtained from $G$ by adding $b'$ to $G$.
We conclude that $G$ and $G'$ are $F_r$-mates. 
\end{proof}

\begin{proof}[Proof of Theorem \ref{th:delta}(i)]
\ \\
(i): We first note that $Forb(C(r,1)) = \bigcirc_r \cup \bigcirc_r^* \cup \cQ^*$.

\noindent
The graphs in $\bigcirc_r$ have maximum degree $D =2$.
By Proposition \ref{le:degree}, almost all graphs in 
$\mathbb{G}(n, d/n)$  have maximum degree $\geq 3$. 
Hence almost all graphs are not in $\bigcirc_r$.

\noindent
By Proposition \ref{pr:8-2} all graphs in $\bigcirc_r^*$ have a $F_r$-mate.
\\
By  Fact \ref{step:8-2} all graphs in $\cQ^*$ have a $F_r$-mate.

\noindent
We conclude that almost all graphs in $\mathbb{G}(n, d/n)$ have a $F_r$-mate.
\end{proof}

\subsection{Proof of Theorem \ref{th:delta}(ii)}
\label{nse:box}
\label{sse:proof2}


\begin{lemma}
\label{nle:1}
Let $G \in \mathbb{G}(n, 1/2)$
be a graph of order $n =rm$ and $A_1, \ldots, A_m$ with each  $A_i \subset V(G), |A_i|=r$, for $i \in [m]$ be a partition of $V(G)$.
Then w.h.p. there exists $i \ [m]$ such that $G[A_i]$ contains a $K_3$.
\end{lemma}
\begin{proof}
Let $G$ be a random graph $G$ in $\mathbb{G}(n, 1/2)$ and $d$ be a constant.
With high probability,   for all $S \subset V(G)$ of size $n/d$, $G[S]$ contains a triangle. 
For this, it suffices to show that the probability that $G[S]$ is triangle-free for a fixed $S$ is $o(1/2^n)$.
If we fix $S$, then $G[S]$ can be viewed as a random graph in  $G \in \mathbb{G}({n}{d}, 1/2)$.

We now partition the vertices $S$ into three sets 
$$
A = \{a_0,\ldots,a_{\frac{n}{3m–1}}\}, B = \{b_0,\ldots,b_{\frac{n}{3m–1}}\}, C = \{c_0,\ldots,c_{\frac{n}{3m–1}}\}.
$$
Consider the collection of triangles 
$$
\{ {a_i, b_j, c_k} : i + j + k = 0 \pmod{\frac{n}{3m}}\}.
$$
This is a collection of $(\frac{n}{3m})^2$ triangles which are edge disjoint: 
for example, any triangle containing the edge $(a_i, b_j)$ must be ${a_i, b_j, c_{–i–j}}$.
The probability that none of these triangles appear in the random graph is $(7/8)^{(n/3m)^2} = o(1/2^n)$, as needed.
\end{proof}

\begin{theorem}
\label{th:box}
In $\mathbb{G}(n, 1/2)$ for almost all graphs $G$ the polynomial $F_r(G;k)$ for $r \geq 4$ vanishes.
\end{theorem}
\begin{proof}
If $G$ has order $n \neq 0 \pmod{r}$ this is 
Proposition \ref{pr:8-1}(v).
If $G$ has order $n = 0 \pmod{r}$, we use Lemma \ref{nle:1}.
If  $F_r(G;k)$
does not vanish, there exists a partition $A_i \subset V(G), |A_i|=r$, for $i \in [m]$ where every
block $A_i$ induces a graph $G[A_i] \cong C_r$. By Lemma \ref{nle:1} this is almost never the case.
\end{proof}

\begin{corollary}
$F_r(G;k)$
with $r \geq4$ has infinitely many unique graphs and is weakly distinguishing. 
In particular, almost all graphs are $F_r$-mates with
\\
$F_r(G;k) =0$ for all $k$. 
Hence, there are infinitely many Harary polynomials which are weakly distinguishing and have infinitely many unique graphs.
\end{corollary}
\begin{proof}
Take $F_r(G;k)$ for $r \geq 4$.
\end{proof}

In Section \ref{se:generating} we briefly looked at graph polynomials $F$ as generating functions and their $F$-mates.
As noted before, Noy's Theorem \ref{pr:noy} states that the independenc polynomial $\rIND(G;k)$ is weakly distinguishing.

For a graph property $\cP$ let $c_{\cP}(G, i)$ be the number of induced subgraphs of $G \in \cP$ of order $i$.
We look at graph polynomials of the form
$$
F_{\cP}(G; x) = \sum_i  c_{\cP}(G; i) x^i.
$$
They are generating functions for the number of induced subgraphs in $\cP$.

\begin{problem}
Find infinitely many graph polynomials of the form of generating functions
$ F_{\cP}(G;k)$ 
which are weakly distinguishing and have infinitely many \\
$F_{\cP}$-unique graphs.
\end{problem}
\ifnewpage
\newpage
\else \fi 

\section{Conclusions and further research}
\label{se:conclu}
\ifmargin
\marginpar{\red File: 
\\
CM-conclu.tex
\\
Last updated:
\\
22.12.25}
\else \fi 

In Section \ref{se:chromatic} I have set the stage for the study of the distinctive power of Harary polynomials.
The various open questions and conjectures can be seen is research program for future work, and I hope that younger researchers
will pick up the thread and work on these problems.
In the following sections I have presented partial results which serve to justify the study of Harary polynomials and their
distinctive power. 
In Section \ref{se:chromatic} I summarized the main results from \cite{HMR}.
The main new results in this paper are
Theorem \ref{th:hereditary},
Theorem \ref{th:gen-hereditary},
Theorem \ref{th:dp-1},
Theorem \ref{th:dp-2} and
Theorem \ref{th:delta} together with Theorem \ref{th:box}.

\subsection*{Acknowledgments}
This paper has a long history. Harary polynomials where introduced by the author already in \cite{makowsky2006polynomial}
and discussed in various graduate courses and later papers  of the author since then.
Vsevolod Rakita's PhD Thesis, \cite{Rakita-PhD},  was dedicated to Harary polynomials and its main results were published as \cite{HMR}.
Preliminary results of this paper were presented at the Workshop on Uniqueness and Discernement in Graph Polynomials in October 2023
at the MATRIX Institute and at the Minisymposium on Graph Polynomials at CanaDAM 2025.
I would like to thank the participants of the Workshop and the Minisymposium for inspiring discussions and
my co-authors O. Herscovici and V. Rakita of our previous paper on Harary Polynomials \cite{HMR} for their collaborative efforts
and contributions.
Finally, I would like to thank 
Orli Herscovici,
Vsevolod Rakita, Elena Ravve, and especially Peter Tittmann, for reading and commenting earlier versions of the paper.
I would also like to thank Yuval Filmus for sharing his expertise on random graphs,
and to the anonymous referees for their valuable suggestions.
\ifnewpage
\newpage
\else \fi 

\newpage
\bibliographystyle{amsplain}
\providecommand{\bysame}{\leavevmode\hbox to3em{\hrulefill}\thinspace}
\providecommand{\MR}{\relax\ifhmode\unskip\space\fi MR }
\providecommand{\MRhref}[2]{%
  \href{http://www.ams.org/mathscinet-getitem?mr=#1}{#2}
}
\providecommand{\href}[2]{#2}

\end{document}